\documentclass[conference, 10pt]{IEEEtran}
\usepackage[style=ieee, backend=biber, isbn=false]{biblatex}

\usepackage{import}

\usepackage[utf8]{luainputenc}

\usepackage[T1]{fontenc}

\usepackage{graphicx}
\usepackage[dvipsnames]{xcolor}

\usepackage[english]{babel}
\addto\captionsenglish{\renewcommand{\figurename}{Fig.}}
\addto\captionsenglish{\renewcommand{\tablename}{Tab.}}
\usepackage{csquotes}

\usepackage{newtxtext}
\usepackage{amsthm}
\usepackage[slantedGreek]{newtxmath}
\usepackage[OMLmathsfit]{isomath}
\DeclareMathAlphabet{\mathbfsf}{\encodingdefault}{\sfdefault}{bx}{n}
\usepackage{bm}
\usepackage{envmath}
\usepackage{mathtools}
\usepackage{commath}
\usepackage{siunitx}

\usepackage[caption=false,font=footnotesize]{subfig}

\usepackage{booktabs}
\usepackage{footmisc}  

\usepackage{url}

\theoremstyle{definition}

\theoremstyle{plain}

\theoremstyle{remark}

\usepackage{lineno}
\modulolinenumbers[5]
\usepackage{umoline}

\usepackage{pgfplots}
\usepackage{pgfplotstable}
\pgfplotsset{compat=newest}
\pgfplotsset{plot coordinates/math parser=false}
\newlength\figureheight
\newlength\figurewidth
\pgfplotsset{every axis plot/.append style={line width=1.5pt},
    legend style={font=\footnotesize, 
        text height=1.0ex,
        draw=black,
        fill=white,
        legend cell align=left}}

\usepackage{hyperref} 
\usepackage[english]{cleveref}

\Crefname{defn}{definition}{definitions}
\Crefname{defn}{Definition}{Definitions}

\Crefname{asm}{assumption}{assumptions}
\Crefname{asm}{Assumption}{Assumptions}

\crefname{lem}{lemma}{lemmas} 
\Crefname{lem}{Lemma}{Lemmas}

\crefname{prop}{proposition}{propositions} 
\Crefname{prop}{Proposition}{Propositions}

\crefname{thm}{theorem}{theorms} 
\Crefname{thm}{Theorem}{Theorms}

\crefname{cor}{corollary}{corollaries}
\Crefname{cor}{Corollary}{Corollaries}
\newcounter{subequation}
\newlength\mtabskip\mtabskip=-1.25cm

\def\mtabLong{long}

\newcommand{\mr}{\mathrm}

\newcommand{\veg}[1]{\bm{#1}}     
\newcommand{\mat}[1]{\mathsfbfit{#1}} 
\renewcommand{\vec}[1]{\mathsfbfit{#1}} 
\newcommand{\op}[1]{\mathcal{#1}} 
\newcommand{\vecop}[1]{\bm{\mathcal{#1}}} 

\newcommand{\dd}{\mathrm{d}}  

\newcommand{\T}{\mr{T}}

\newcommand\restr[2]{{
        \left.\kern-\nulldelimiterspace 
        #1 
        \vphantom{|} 
        \right|_{#2} 
}}

\newcommand\rst[3]{{
        \left.\kern-\nulldelimiterspace 
        #1 
        \vphantom{|} 
        \right|_{#2}^{#3} 
}}

\usepackage{acro}

\DeclareAcronym{DG}
{
    short = DG ,
    long = discontinuous Galerkin
}

\DeclareAcronym{ACA}
{
    short = ACA ,
    long = adaptive cross approximation
}

\DeclareAcronym{EFIE}
{
    short =  EFIE ,
    long = electric field integral equation
}

\DeclareAcronym{MFIE}
{
    short =  MFIE ,
    long = magnetic field integral equation
}

\DeclareAcronym{CFIE}
{
    short =  CFIE ,
    long = combined field integral equation
}

\DeclareAcronym{MUIE}
{
    short =  MUIE ,
    long = Müller integral equation
}

\DeclareAcronym{PMCHWT}
{
    short =  PMCHWT ,
    long = Poggio-Miller-Chang-Harrington-Wu-Tsai integral equation
}

\DeclareAcronym{SPD}
{
    short =  SPD ,
    long = {symmetric, positive definite}
}

\DeclareAcronym{SPSD}
{
    short =  SPD ,
    long = {symmetric, positive semi-definite}
}

\DeclareAcronym{PEC}
{
    short =  PEC ,
    long = perfectly electrically conducting
}

\DeclareAcronym{RWG}
{
    short = RWG ,
    long = Rao-Wilton-Glisson
} 

\DeclareAcronym{BC}
{
    short = BC ,
    long = Buffa-Christiansen
}

\DeclareAcronym{SVD}
{
    short = SVD ,
    long = singular value decomposition
}

\DeclareAcronym{CG}
{
    short = CG ,
    long = conjugate gradient
} 

\DeclareAcronym{PCG}
{
    short = PCG ,
    long = preconditioned conjugate gradient
} 

\DeclareAcronym{CGS}
{
    short = CGS ,
    long = conjugate gradient squared
}

\DeclareAcronym{CMP}
{
    short = CMP ,
    long = Calderón multiplicative preconditioner
} 

\DeclareAcronym{RFCMP}
{
    short = RF-CMP ,
    long = refinement-free Calderón multiplicative preconditioner
} 

\DeclareAcronym{HPD}
{
    short = HPD ,
    long = {Hermitian, positive definite}
} 

\DeclareAcronym{RHS}
{
    short = RHS ,
    long = {right-hand side}
}

\DeclareAcronym{PW}
{
    short = PW ,
    long = {plane wave}
} 

\DeclareAcronym{HD}
{
    short = HD ,
    long = {Hertzian dipole}
} 

\DeclareAcronym{FF}
{
    short = FF ,
    long = {far-field}
} 

\DeclareAcronym{NF}
{
    short = NF ,
    long = {near-field}
}  

\newcolumntype {n}{c}
\newcolumntype {N}{>{\small}c}
\newcolumntype {L}{>{\small}l}
\newcolumntype {F}{>{\footnotesize}c}
\newcolumntype {v}[1]{>{\raggedright \hspace {0pt}} p {#1}}
\newcolumntype {V}[1]{>{\small \raggedright \hspace {0pt}} p {#1}}
\newcolumntype{d}[1]{>{\DC@{.}{.}{#1}}c<{\DC@end}}

\newcolumntype{R}[1]{%
    >{\begin{turn}{90}\begin{minipage}{#1}\small\raggedright\hspace{0pt}}l%
            <{\end{minipage}\end{turn}}%
}

\NewDocumentCommand{\TA}{o}{
    \IfNoValueTF {#1} {%
        \vecop T_{\kern-2pt\mr{A}}
    }
    {
        \vecop T_{\kern-2pt\mr{A},#1}
    }
}

\NewDocumentCommand{\TPhi}{o}{
    \IfNoValueTF {#1} {%
        \vecop T_{\kern-2pt\Phiup}
    }
    {
        \vecop T_{\kern-2pt\Phiup,#1}
    }
}

\NewDocumentCommand{\matTA}{o}{
    \IfNoValueTF {#1} {%
        \mat T_\mr{A}   
        }
    {
        \mat T_{\mr{A},#1}
    }
}

\NewDocumentCommand{\matTPhi}{o}{
    \IfNoValueTF {#1} {%
        \mat T_\Phiup   
        }
    {
        \mat T_{\Phiup,#1}
    }
}

\NewDocumentCommand{\MSL}{o}{
    \IfNoValueTF {#1} {%
        \veg \Psi_\mr{SL}
        }
    {
        \veg \Psi_{\mr{SL},#1}
    }
}

\NewDocumentCommand{\MDL}{o}{
    \IfNoValueTF {#1} {%
        \veg \Psi_\mr{DL}
        }
    {
        \veg \Psi_{\mr{DL},#1}
    }
}

\NewDocumentCommand{\PA}{o}{
    \IfNoValueTF {#1} {%
        \veg \Psi_\mr{A}
        }
    {
        \veg \Psi_{\mr{A},#1}
    }
}

\NewDocumentCommand{\PPhi}{o}{
    \IfNoValueTF {#1} {%
        \veg \Psi_{\Phiup}
        }
    {
        \veg \Psi_{\Phiup,#1}
    }
}

\usepackage{comment}
\usepackage[belowskip=-4pt,aboveskip=0pt,font=footnotesize]{caption}

\newcommand\LapLm[1][]{\mat{\Delta}_{#1}}
\newcommand{\LapLim}[1][]{\mat{\Delta}^{-1}_{#1}}
\newcommand\LapSm[1][]{{\tilde{\mat{\Delta}}_{#1}}}

\newcommand\patch[1][]{{\pi_{#1}}}

\newcommand\pyr[1][]{{\lambda_{#1}}}
\newcommand\pyrD[1][]{{\tilde{\lambda}_{#1}}}

\newcommand\GmixR[1][]{\mat{G}_{\pyrD[#1]\patch[#1]}}
\newcommand\LapSP[1][]{{\tilde{\mat{\Delta}}_{\pi#1}}}

\def\Nm{\mat{N}}
\def\Dm{\mat{D}}
\def\Sm{\mat{S}}
\def\Dxm{\mat{D}^*}
\pdfmapfile{-mpfonts.map}

\makeatletter

\newcounter{authr}
\newcommand{\authr}[2][]{
   \stepcounter{authr}
   \@namedef{authr@\theauthr}{#2}
   \@namedef{authrlabel@\theauthr}{#1}
}

\newcounter{address}
\newcommand{\address}[2][]{
   \stepcounter{address}
   \@namedef{address@\theaddress}{#2}
   \@namedef{addresslabel@\theaddress}{#1}
}

\newcommand{\alsep}{and}

\def\newmaketitle{\par%
  \begingroup%
  \normalfont%
  \def\thefootnote{}
  \def\footnotemark{}
  \let\@makefnmark\relax
  \footnotesize
  \footnotesep 0.7\baselineskip
  \normalsize%
  \twocolumn[\thenewmaketitle\@IEEEaftertitletext]%
  \if@IEEEusingpubid
     \enlargethispage{-\@IEEEpubidpullup}%
  \fi
  \endgroup
  \setcounter{footnote}{0}\let\maketitle\relax\let\@maketitle\relax
  \gdef\@thanks{}%
  \let\thanks\relax}

\def\thenewmaketitle{
  \newpage
  \begin{center}%
    \vskip0.2em{\Huge\@IEEEcompsoconly{\sffamily}\@IEEEcompsocconfonly{\normalfont\normalsize\vskip 2\@IEEEnormalsizeunitybaselineskip
   \bfseries\large}\@title\par}\vskip1.0em\par%
    \vspace{1ex}
    \newcounter{c@authr}
    \newcounter{c@tmp}
    \ifthenelse{\value{authr}=2}{%
      \newcommand{\liand}{ and }}{%
      \newcommand{\liand}{, and }}
    \ifthenelse{\value{address}<2}{%
      \@nameuse{authr@1}%
      \stepcounter{c@authr}%
      \whiledo{\value{c@authr}<\value{authr}}{%
        \setcounter{c@tmp}{\value{authr}}%
        \addtocounter{c@tmp}{-\value{c@authr}}%
        \ifthenelse{\value{c@tmp}=1}{%
          \renewcommand{\alsep}{\liand}}{\renewcommand{\alsep}{, }}%
        \stepcounter{c@authr}\alsep \@nameuse{authr@\thec@authr}}\\%
    }
    {
      \@nameuse{authr@1}${}^{(\ref{\@nameuse{authrlabel@1}})}$%
      \stepcounter{c@authr}%
      \whiledo{\value{c@authr}<\value{authr}}{%
      \setcounter{c@tmp}{\value{authr}}%
      \addtocounter{c@tmp}{-\value{c@authr}}%
      \ifthenelse{\value{c@tmp}=1}{%
        \renewcommand{\alsep}{\liand}}{\renewcommand{\alsep}{, }}%
      \stepcounter{c@authr}\alsep \@nameuse{authr@\thec@authr}%
        ${}^{(\ref{\@nameuse{authrlabel@\thec@authr}})}$%
      }
    }
    \vspace{0.2ex}

    \ifthenelse{\value{address}>0}{%
      \ifthenelse{\value{address}=1}{
        {\@nameuse{address@1}}
      }
      {
        \newcounter{c@address}

        \begin{center}
        \whiledo{\value{c@address}<\value{address}}
        {
          \refstepcounter{c@address}
            ${}^{(\thec@address)}$\,%
              \label{\@nameuse{addresslabel@\thec@address}}%
              \@nameuse{address@\thec@address}\\ %
        }
        \end{center}
      } 
    }
    {
      \relax
    }
  \end{center}
}

\makeatother

\title{A New Refinement-Free  Preconditioner for the Symmetric Formulation in Electroencephalography
}

\authr[org1]{Viviana Giunzioni}
\authr[org2]{John E. Ortiz G.}
\authr[org3]{Adrien Merlini}
\authr[org4]{Simon B. Adrian}
\authr[org1]{Francesco P. Andriulli}

\address[org1]{Department of Electronics and Telecommunications, Politecnico di Torino, 10129 Turin, Italy}
\address[org2]{Electronic Engineering Department, Universidad de Nariño, Pasto 520002, Colombia}
\address[org3]{Microwaves Department, IMT Atlantique, 29238 Brest, France}
\address[org4]{Fakultät für Informatik und Elektrotechnik, Universität Rostock, 18051 Rostock, Germany}

\begin{document}

\newmaketitle

\begin{abstract}
Widely employed for the accurate solution of the electroencephalography forward problem, the symmetric formulation gives rise to a first kind, ill-conditioned operator ill-suited for complex modelling scenarios.
This work presents a novel preconditioning strategy based on an accurate spectral analysis of the operators involved which, differently from other Calder\'{o}n-based approaches, does not necessitate the barycentric refinement of the primal mesh (i.e., no dual matrix is required).
The discretization of the new formulation gives rise to a well-conditioned, symmetric, positive-definite system matrix, which can be efficiently solved via fast iterative techniques.
Numerical results for both canonical and realistic head models validate the effectiveness of the proposed formulation.
\end{abstract}

\section{Introduction}

Electroencephalographic (EEG) data, coupled with a reliable source localization algorithm, provide a valuable tool for the non-invasive mapping of the brain activity \cite{darbas2019review}.
A crucial application of this technology is the pre-surgical characterization of epileptic seizures, plaguing approximately $1\%$ of the world population \cite{darbas2019review}.
The localization of the active neural sources inside the brain from a measure of the potential at the scalp---referred to as the inverse EEG problem---relies on iterative solutions of the forward EEG problem, which, as a consequence, directly affects the accuracy and resolution time \cite{darbas2019review} of the method.
This consideration sheds light on the paramount need for solving accurately \textit{and} efficiently the forward EEG problem. Among the many integral equation strategies available, the symmetric formulation \cite{kybic2005common} has become particularly popular because of the higher accuracy it yields.
Unfortunately, it suffers from severe numerical instabilities when the spatial resolution becomes finer and when the conductivity contrast between adjacent layers increases, thus compromising its application to realistic scenarios.

Recently, a novel refinement-free preconditioning strategy based on pseudo-differential operator theory has been successfully applied to the electric field integral equation \cite{adrian2019refinement}. In a similar philosophy, we developed a new preconditioning scheme for the symmetric formulation, which, differently from standard Calder\'{o}n-based approaches \cite{pillain2018calderon}, does not require the evaluation of the electromagnetic operators over dual functions defined on the barycentric refinement of the primal mesh, thus considerably reducing the building time of the boundary element method (BEM) matrices. Moreover, the symmetric, positive definiteness of the formulation makes it amenable to fast iterative solvers.

\section{Background and Notation}

The forward EEG problem consists of finding the potential $V$ at the positions of the electrodes generated by a known configuration of current sources $\bm{j}$. Given the head conductivity profile $\sigma$, the problem is modelled by Poisson's equation $\nabla \cdot (\sigma \nabla V) = \nabla \cdot \bm{j}$ \cite{darbas2019review}.
Let $\Omega = \cup_{i=1}^{N}\Omega_i$ be a nested domain with boundaries $\Gamma_i = \overline{\Omega}_i\cap\overline{\Omega}_{i+1}$ characterized by the outgoing normal fields $\veg{n}_i$. By assuming the conductivity $\sigma_i$ constant in each $\Omega_i$, the application of the representation theorem yields a system of $2N$ equations \cite{kybic2005common}
\begin{align}
&(\uppartial_{\veg n} v_{i+1})_{\Gamma_i}-(\uppartial_{\veg n} v_{i})_{\Gamma_i} = -\op{D}^*_{i,i-1}p_{i-1}+2\op{D}^*_{ii}p_i-\op{D}^*_{i,i+1}p_{i+1}
\nonumber \\
&+\sigma_i \op{N}_{i,i-1}V_{i-1}-(\sigma_i+\sigma_{i+1})\op{N}_{ii}V_i+\sigma_{i+1}\op{N}_{i,i+1}V_{i+1}
\\
&\sigma_{i+1}^{-1}(v_{i+1})_{\Gamma_i}-\sigma_{i}^{-1}(v_{i})_{\Gamma_i} = \mathcal{D}_{i,i-1}V_{i-1}-2\mathcal{D}_{ii}V_i+\mathcal{D}_{i,i+1}V_{i+1}
\nonumber \\
&-\sigma_{i}^{-1}\mathcal{S}_{i,i-1}p_{i-1}+(\sigma_i^{-1}+\sigma_{i+1}^{-1})\mathcal{S}_{ii}p_i-\sigma_{i+1}^{-1}\mathcal{S}_{i,i+1}p_{i+1}
,
\end{align}
where $V_i = (V)_{\Gamma_i}$, $p_i = \sigma_i[\veg{n}_i\cdot \nabla V]_{\Gamma_i}$, and
$v_i = f_i*G$, with $f_i$ equal to $\nabla \cdot \bm{j}$ in $\Omega_i$ and to $0$ elsewhere and with the free-space Green's function $G(\veg{r}, \veg{r}') = 1/(4\uppi |\veg{r}-\veg{r}'|)$.
The operators involved are the single-layer $(\mathcal{S}\psi)(\veg{r}) = \int_{\uppartial \Omega}G(\veg{r}, \bm{r}') \psi(\veg{r}') \dd S(\veg{r}')$, the double layer $\mathcal{D}\phi(\veg{r}) = \mathrm{p.v.} \int_{\uppartial \Omega}\uppartial_{\veg{n}'} G(\veg{r}, \veg{r}') \phi(\veg{r}') \dd S(\veg{r}')$, the adjoint double layer $\mathcal{D}^*\psi(\veg{r}) = \mathrm{p.v.} \int_{\uppartial \Omega}\uppartial_{\veg{n}} G(\veg{r}, \veg{r}') \psi(\veg{r}') \dd S(\veg{r}')$ and the hypersingular $\mathcal{N}\phi(\veg{r}) = \int_{\uppartial \Omega}\uppartial_{\veg{n}}\uppartial_{\veg{n}'}G(\veg{r}, \veg{r}') \phi(\veg{r}')\dd S(\veg{r}')$ operators; the subscript $_{ij}$ denotes the restriction of the operator to $\veg{r}\in\Gamma_i$ and $\veg{r}'\in\Gamma_j$.
The above system can be solved numerically on a mesh with $N_\mr{c}$ cells and $N_\mr{v}$ vertices, over which we define a set of piecewise constant, $\{\patch[n]\}_{n=1}^{N_\mr{c}}$, and piecewise linear, $\{\pyr[n]\}_{n=1}^{N_\mr{v}}$, basis functions as in \cite{adrian2019refinement} in equations (2), (3).
To proceed with the numerical solution of the system, the unknowns are expanded as $V_i=\sum_{n=1}^{N_\mr{v}}{v}_{i,n} \pyr[n]$, $p_i=\sum_{n=1}^{N_\mr{c}}{p}_{i,n} \patch[n]$ and the resulting equations are tested with $\pyr$ and $\patch$ functions.
For the case $N=1$, the system finally reads $\mat{Z}[ \vec{v}^\T \,\, \vec{p}^\T ]^\T =[ \vec{b}^\T \,\, \vec{c}^\T ]^\T$, with $\mat{Z} =\big[\mat{N} \,\, \mat{D}^*; \mat{D} \,\, \mat{S}\big]$,
in which the blocks are
$\sbr{\mat{N}}_{mn}= (\sigma_1+\sigma_2)\left( \pyr_m, \op{N}_{11}\pyr_n \right)_{L^2(\Gamma_1)}$,
$\sbr{\mat{D}^*}_{mn} = -2\left( \pyr_m, \op{D}^*_{11}\patch_n \right)_{L^2(\Gamma_1)}$,
$\sbr{\mat{D}}_{mn} = -2\left( \patch_m, \op{D}_{11}\pyr_n \right)_{L^2(\Gamma_1)}$, and
$\sbr{\mat{S}}_{mn} = (\sigma_1^{-1}+\sigma_2^{-1})\left( \patch_m, \op{S}_{11}\patch_n \right)_{L^2(\Gamma_1)}$,
and where $\sbr{\vec{v}}_m = {v}_{1,m}$, $\sbr{\vec{p}}_m = {p}_{1,m}$, $\sbr{\vec{b}}_m = -(\pyr[m], \uppartial_{\veg n} v_1)_{L^2(\Gamma_1)}$ and $\sbr{\vec{c}}_m = -\sigma_1^{-1}(\patch[m], v_1)_{L^2(\Gamma_1)}$.

\section{Our New Formulation}
Two different sources of ill-conditioning plague the symmetric formulation. On the one hand, we observe numerical instabilities related to the high conductivity contrast between adjacent regions, necessarily present in realistic head models; these are cured by applying a proper rescaling with respect to the conductivity \cite{pillain2018calderon}. On the other hand, the symmetric formulation suffers from a severe dense-discretization breakdown caused by the opposed asymptotic behaviours of the spectra of the singular and of the hypersingular operators.
The remedy we propose here is a preconditioning based on the regularization of an operator of pseudo-differential order $n$ via the multiplication by other operators of orders summing up to $-n$. In Calder\'{o}n approaches, this product chain is limited to only two elements, for example $\op{N}$ is regularized by the left multiplication by $\op{S}$ and vice versa. Instead, significant advantages follow from a three-operator multiplication, where the first and the last elements are the same. Consider for example the block $\mat{N} \mat{X} \mat{N}$, where $\mat{X}$ is the discretization of the unknown operator $\op{X}$: provided that the hypersingular operator has pseudo-differential order $1$, the above will be well-conditioned for any $\op{X}$ of order $-2$. Moreover, the symmetric positive semi-definiteness (SPsD) of $\mat{X}$ would also guarantee the same property for $\mat{N} \mat{X}  \mat{N}$. Following this philosophy, we designed our new formulation that reads, for the sake of brevity limited here to the one-compartment case only,
\begin{equation}
    \mat M \, \mat{Z} \,\mat{P} \,\mat{Z} \mat{M} \,\,\vec{y} =  \mat{M}\, \mat{Z}\, \mat{P} \, \begin{bmatrix} \vec{b}^\T & \vec{c}^\T \end{bmatrix}^\T.
    \label{eqn:completeFormulation}
\end{equation}
Once the system in \eqref{eqn:completeFormulation} is properly deflated, the solution of the original problem is retrieved as $[ \vec{v}^\T \,\, \vec{p}^\T ]^\T=\mat{M} \vec{y}$.
The matrices $\mat{M}$ and $\mat{P}$ are defined as $\mat{M} = \text{diag}(\mat{G}^{-1/2}_{\lambda\lambda}, \mat{G}^{-1/2}_{\pi\pi})$, $\mat{P}=\text{diag}(\LapLim[],\LapSP[])$,
where $\mat{G}_{\lambda\lambda}$ and $\mat{G}_{\pi\pi}$ are the pyramid and patch Gram matrices and $\LapSP[] = \GmixR[]^{-1}\LapSm[]\GmixR[]^{-1}$, with $\GmixR[]$ the Gram matrix linking the dual pyramid and the patch functional spaces.
The blocks $\LapLm[]$ and $\LapSm[]$ discretize the Laplace-Beltrami operator $\Deltaup$ with primal and dual piecewise linear functions respectively. 
It can be proven that the result arising from
\begin{equation}
    \mat Z \mat P \mat Z = 
     \bigg[
    \begin{array}{cc}
        \Nm\LapLim\Nm + \Dxm\LapSP\Dm & \Nm\LapLim\Dxm + \Dxm\LapSP\Sm \\
        \Dm\LapLim\Nm + \Sm\LapSP\Dm & \Dm\LapLim\Dxm + \Sm\LapSP\Sm
    \end{array}
  \bigg]
\end{equation}
is symmetric positive semi-definite and spectrally equivalent to a four identity-block matrix, plus a compact perturbation. 
This is compatible with the fact that, by assigning an order to each matrix corresponding to the pseudo-differential order of the underlying operator and by summing the orders of the concatenated blocks---recalling that $\op{N}$, $\op{S}$, $\op{D}$, $\op{D}^*$, ${\Deltaup}$ have orders $1$, $-1$, $-1$, $-1$ and $2$ respectively---all dominant blocks' orders sum to zero.
A careful analysis then obtained by evaluating analytically the eigenvalues of the blocks, ensures that the principal part of the operator is not deleted and  an asymptotic bound on the condition number can be obtained, ensuring the well-posedness of the formulation. 

\section{Numerical Results}
\begin{figure}
\centerline{\includegraphics[width=0.82\columnwidth]{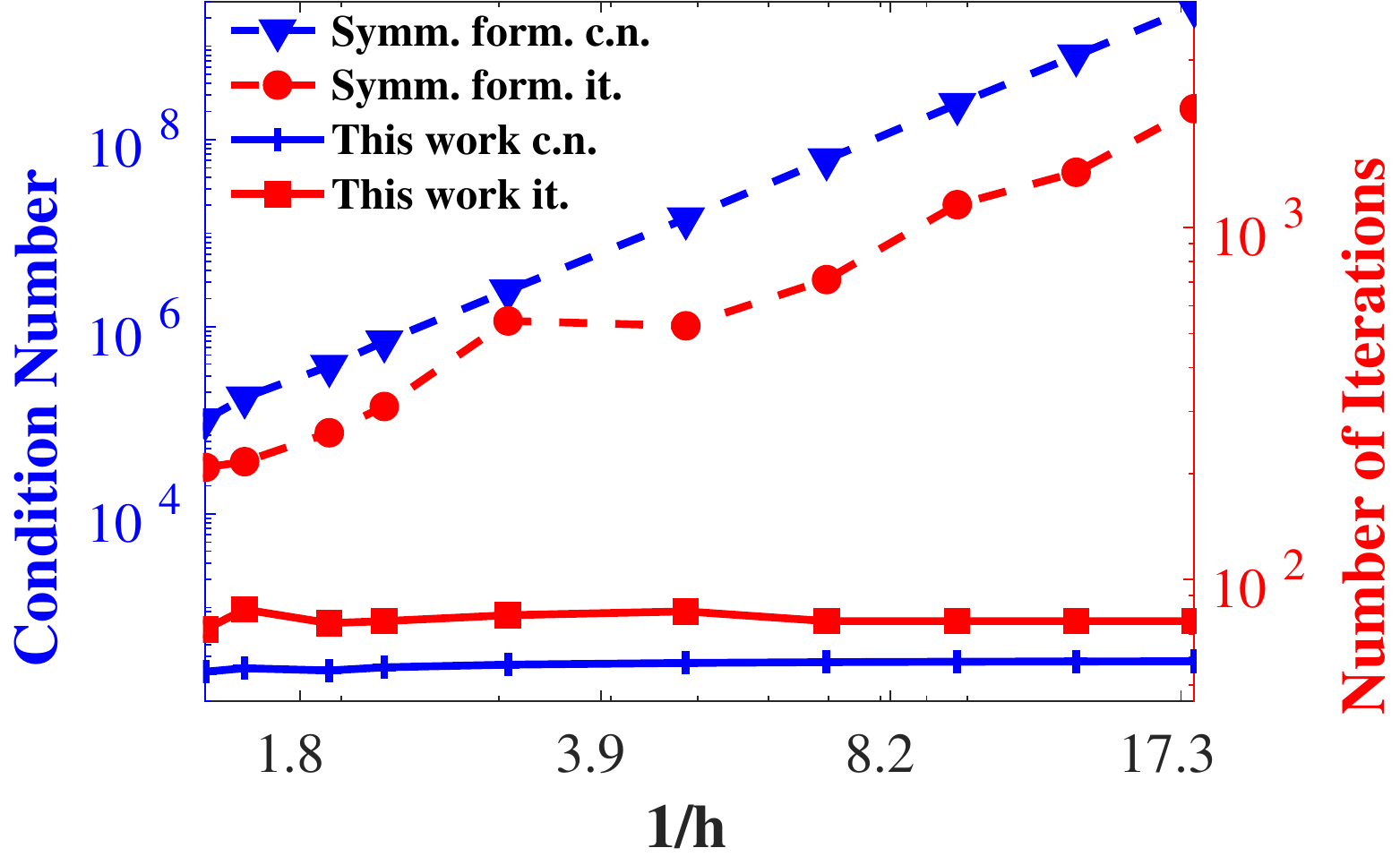}}
\vspace*{0mm}
\caption{Condition number and number of iterations as a function of the inverse average edge length $1/h$.}
\vspace*{-2mm}
\label{fig:cn}
\end{figure}

The stability of the new formulation has been tested on three concentric spheres with radii ($0.8$, $0.9$, $1$) in normalized units, modelling the brain, the skull, and the skin, for which the analytic solution for the potential on the outermost layer is available \cite{darbas2019review}. As is often the case in the literature \cite{darbas2019review}, the conductivities of the three regions have been set to $1{:}1/80{:}1$ and the system has been excited through a current dipole 
source at $0.6$ from the center and radially directed. Clearly, the results in Fig.~\ref{fig:cn} show that both the condition number and the number of iterations remain constant for the new formulation, differently from the standard symmetric formulation case.

\begin{figure}
\subfloat[\label{subfig-1:headA}]{%
  \includegraphics[width=0.53\columnwidth]{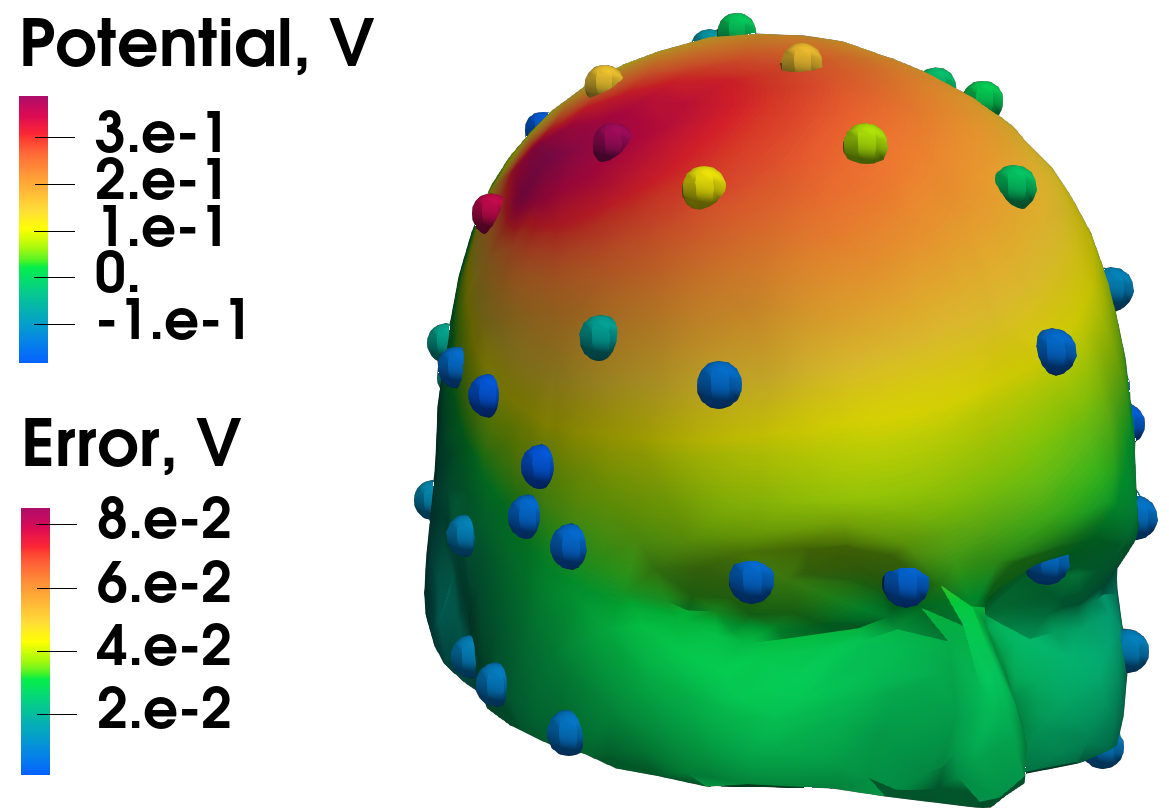}
}
\hfill
\subfloat[\label{subfig-2:headB}]{%
  \includegraphics[width=0.334\columnwidth]{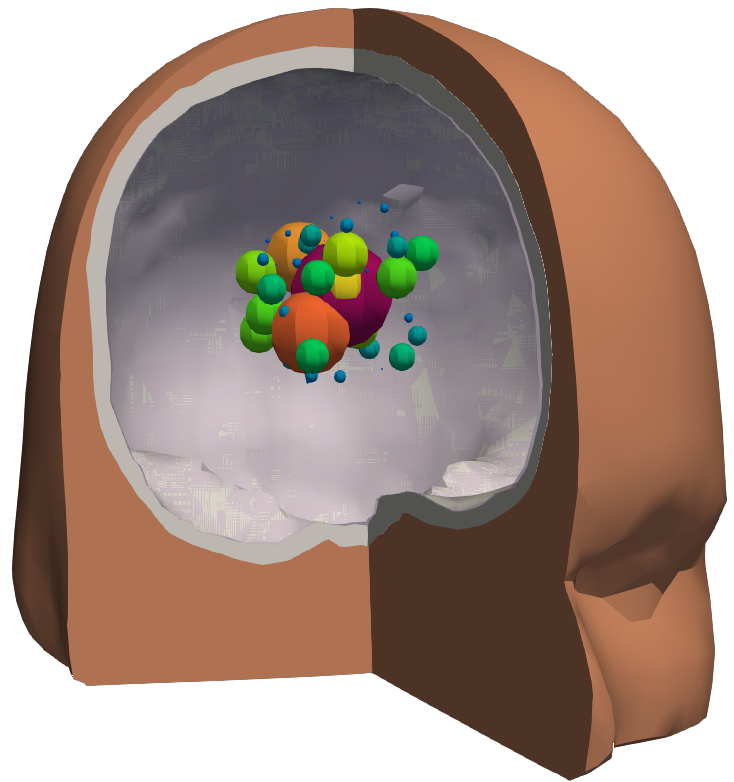}
}
\vspace*{0.5mm}
\caption{(a) Scalp potential distribution and error at the electrodes with respect to the symmetric formulation solution. (b) Reconstructed epileptogenic source.}
\end{figure}

We have also applied our formulation on a realistic head model extracted from magnetic resonance imaging data, to obtain the scalp potential distribution in  Fig.~\ref{subfig-1:headA}. Given the good agreement with the reference solution, our formulation can be effectively used to solve an inverse EEG problem, such as the localization of epileptogenic zones, as shown in Fig.~\ref{subfig-2:headB}.

%
%
%
%
\section*{Acknowledgment}\label{sec:ack}
This work was supported in part by the European Research Council (ERC) through the European Union’s Horizon 2020 Research and Innovation Programme under Grant 724846 (Project 321) and in part  by the Italian Ministry of University and Research within the Program FARE, CELER, under Grant R187PMFXA4

\vspace{-0.15cm}


{\small \printbibliography}
\end{document}